\documentclass[a4paper, final, 12pt]{article}
\usepackage{amsmath, amssymb, latexsym, amscd, amsthm,amsfonts,amstext}
\usepackage[mathscr]{eucal}
\usepackage{graphicx}
\usepackage{subfig}
\usepackage{float}
\usepackage{color}
\usepackage{hyperref}
\usepackage[utf8]{inputenc}
\usepackage[english]{babel}

\numberwithin{equation}{section}
\input{tcilatex}
\begin{document}

\title{H\"{o}lder Stability and Uniqueness for The Mean Field Games System
via Carleman Estimates\thanks{\textbf{Funding}. The work of J. Li was
partially supported by the NSF of China No. 11971221, Guangdong NSF Major
Fund No. 2021ZDZX1001, the Shenzhen Sci-Tech Fund No. RCJC20200714114556020,
JCYJ20200109115422828 and JCYJ20190809150413261, National Center for Applied
Mathematics Shenzhen, and SUSTech International Center for Mathematics. The
work of H. Liu was supported by the Hong Kong RGC General Research Funds
(projects 12302919, 12301420 and 11300821) and the France-Hong Kong ANR/RGC
Joint Research Grant, A-CityU203/19.}}
\author{ Michael V. Klibanov \thanks{
Department of Mathematics and Statistics, University of North Carolina at
Charlotte, Charlotte, NC, 28223, USA, mklibanv@uncc.edu}, \and Jingzhi Li 
\thanks{
Department of Mathematics \& National Center for Applied Mathematics
Shenzhen \& SUSTech International Center for Mathematics, Southern
University of Science and Technology, Shenzhen 518055, P.~R.~China,
li.jz@sustech.edu.cn}, \and Hongyu Liu \thanks{
Department of Mathematics, City University of Hong Kong, Kowloon, Hong Kong
SAR, P.R. China, hongyliu@cityu.edu.hk}}
\date{}
\maketitle

\begin{abstract}
We are concerned with the mathematical study of the Mean Field Games system
(MFGS). In the conventional setup, the MFGS is a system of two coupled
nonlinear parabolic PDEs of the second order in a backward-forward manner,
namely one terminal and one initial conditions are prescribed respectively
for the value function and the population density. In this paper, we show
that uniqueness of solutions to the MFGS can be guaranteed if, among all
four possible terminal and initial conditions, either only two terminal or
only two initial conditions are given. In both cases H\"{o}lder stability
estimates are proven. This means that the accuracies of the solutions are
estimated in terms of the given data. Moreover, these estimates readily
imply uniqueness of corresponding problems for the MFGS. The main
mathematical apparatus to establish those results is two new Carleman
estimates, which may find application in other contexts associated with
coupled parabolic PDEs.
\end{abstract}


ill-posed and inverse problems, Carleman estimates


\section{Introduction}

\label{sec:1}

The mean field games (MFG) theory was first developed in the seminal works
of Lasry and Lions \cite{LL1,LL2,LL3} as well as of Huang, Caines and Malham%
\'{e} \cite{Huang1,Huang2}. This theory studies the behavior of infinitely
many agents, who are trying to optimize their values. There are many
applications of this theory in \underline{social sciences}. Some examples of
those applications are: finance, economics, pedestrians flocking and
interactions of electrical vehicles, see, e.g. \cite{A,Cou,Kol,LL3,Trusov}.\
We also mention applications of the MFG theory in the fight with corruption 
\cite{KM}, \cite[Preface]{Kol} and the cyber security \cite{KB}.

The methodology of this publication comes from the theory of Ill-Posed and
Inverse Problems. The authors have been heavily involved in this field
throughout their careers, see, e.g. \cite{BK,BukhKlib}, \cite{Klib84}-\cite%
{KTnum}, \cite{LL1,LLLZ1}. The apparatus of Carleman estimates was first
introduced in the field of coefficient inverse problems in the work of
Bukhgeim and Klibanov \cite{BukhKlib}, and it was first introduced in the
MFG theory by Klibanov and Averboukh \cite{MFG1}, also, see two follow up
publications \cite{MFG2,MFG3}.

Accuracy estimates for the solution of the mean field games system (MFGS)
with respect to the input data, which the authors also call
\textquotedblleft stability estimates", were unknown in the MFG theory prior
to \cite{MFG1,MFG2,MFG3}. In the meantime such estimates are quite desirable
ones since the input data for the MFGS are given with errors. Besides, these
estimates imply uniqueness of corresponding problems for the MFGS.

Let $u\left( x,t\right) $ be the value function of a mean field game and let 
$m\left( x,t\right) $ be the function, which describes the density of agents 
\cite{A,LL3}. Here $x\in \mathbb{R}^{n}$ and $t$ are the spatial variable
and the time variable respectively and $t\in \left( 0,T\right) .$ The mean
field games system (MFGS) of the second order is a crucial part of the MFG
theory. The MFGS is a system of two coupled nonlinear parabolic equations
with respect to functions $u\left( x,t\right) $ and $m\left( x,t\right) .$\
A substantial complication of the MFGS is that times are running in two
different directions in those two PDEs. Therefore, the conventional theory
of parabolic PDEs is inapplicable to the MFGS.

The following are four possible terminal and initial conditions for the MFGS:%
\begin{equation}
u\left( x,T\right) ,m\left( x,0\right) ,m\left( x,T\right) ,u\left(
x,0\right) .  \label{1.1}
\end{equation}%
In the conventional setting the following two out of these four functions
are known: 
\begin{equation}
u\left( x,T\right) \text{ and }m\left( x,0\right) .  \label{1.2}
\end{equation}%
In addition, usually functions $u\left( x,t\right) $ and $m\left( x,t\right) 
$ are assumed to be periodic with respect to each component of the vector $x$%
, see, e.g. \cite{A,LL3}.\ However, uniqueness of the solution of the MFGS
is in question then, unless quite restrictive the so-called
\textquotedblleft monotonicity" conditions are not imposed \cite{Bardi}. On
the other hand, if one assumes that either of two vector functions%
\begin{equation}
\left( u\left( x,T\right) ,m\left( x,0\right) ,m\left( x,T\right) \right) ,
\label{1.3}
\end{equation}%
\begin{equation}
\left( u\left( x,T\right) ,m\left( x,0\right) ,u\left( x,0\right) \right)
\label{1.4}
\end{equation}%
\ is known and zero Neumann boundary conditions are imposed on both
functions $u\left( x,t\right) $\ and $m\left( x,t\right) $, then the
Lipschitz stability estimate for either of cases (\ref{1.3}) or (\ref{1.4})
holds along with the uniqueness \cite{MFG1,MFG2}.\ In \cite{MFG3}\ both a H%
\"{o}lder stability estimate and uniqueness are obtained in the case when
Dirichlet and Neumann boundary conditions for functions $u\left( x,t\right) $%
\ and $m\left( x,t\right) $ are known, whereas all functions in (\ref{1.1})
are unknown.

Still, three functions are known in either (\ref{1.3}) or (\ref{1.4}), which
means an over-determination in the data. Hence, the following question is
natural one to be posed: \emph{Can only two out of four terminal and initial
conditions (\ref{1.1})\ provide both a stability estimate and uniqueness for
the MFGS}?\ The goal of this paper is to address this question positively.
More precisely, we demonstrate here that if Neumann boundary conditions are
imposed on both functions $u\left( x,t\right) $ and $m\left( x,t\right) $,
then the replacement of the conventional pair (\ref{1.2}) with either of two
pairs%
\begin{equation}
\left( u\left( x,T\right) ,m\left( x,T\right) \right) \text{ or }\left(
u\left( x,0\right) ,m\left( x,0\right) \right)  \label{1.5}
\end{equation}%
leads to both: a H\"{o}lder stability estimate and uniqueness. In other
word, if the conventional pair (\ref{1.2})\ is replaced by one of two pairs (%
\ref{1.5}), then uniqueness of the solution of the MFGS is restored even for
the non-overdetermined case, and, in addition, H\"{o}lder estimate of the
accuracy of the solution is in place.\ 

Assuming that the function $m\left( x,T\right) $ is known, we actually
assume that we can measure the final distribution of players. Next, solving
the MFGS with the terminal data $\left( u\left( x,T\right) ,m\left(
x,T\right) \right) ,$ we provide a retrospective analysis of the process 
\cite{MFG1}. On the other hand, an approximate knowledge of the\ initial
condition $u\left( x,0\right) $ of the value function can be obtained via
polling of players in the beginning of the process about their ideas about
their value function \cite{MFG2}. At the same time, since polls are
expensive efforts, then it is reasonable to obtain the result of a poll only
once, rather than conducting polls at several moments of time.\ \ \ \ \ \ \
\ \ \ \ \ \ \ \ \ \ \ \ \ \ \ \ \ \ \ \ \ \ \ \ \ \ \ \ \ \ \ \ \ \ \ \ \ \
\ \ \ \ \ \ \ \ \ \ \ \ \ \ \ \ \ \ \ \ \ \ \ \ \ \ \ \ \ \ \ \ \ \ \ \ \ \
\ \ \ \ \ \ \ \ \ \ \ \ \ \ \ \ \ \ \ \ \ \ \ \ \ \ \ \ \ \ \ \ \ \ \ \ \ \
\ \ \ \ \ \ \ \ \ \ \ \ \ \ \ \ \ \ \ \ \ \ \ \ \ \ \ \ \ \ \ \ \ \ \ \ \ \
\ \ \ \ \ \ \ \ \ \ \ \ \ \ \ \ \ \ \ \ \ \ \ \ \ \ \ \ \ \ \ \ \ 

Results of the current publication as well as of \cite{MFG1,MFG2,MFG3} are
about a single measurement case. As to the case of infinitely many
measurements, we refer to two recent results of \cite{Liu1,Liu2}, which
prove uniqueness of the reconstruction of the interaction term of the MFGS.

We rely below on two new Carleman estimates for the MFGS, which were derived
in \cite{MFG1,MFG2}. Carleman estimates are traditionally used for proofs of
stability and uniqueness theorems for ill-posed Cauchy problems for various
PDEs, although only the case of a single PDE is usually considered, unlike
our case of a system of PDEs, see, e.g. see, e.g. \cite%
{KT,Ksurvey,KL,LRS,Yam}. Starting from the originating publication \cite%
{BukhKlib}, Carleman estimates have been actively used for proofs of global
uniqueness and stability results for coefficient inverse problems. Since
this paper is not a survey of publications devoted to the method of \cite%
{BukhKlib}, we refer now only to a few of those and references cited therein 
\cite{BK,Fu,ImYam,Isakov,Klib84,Kl84,Klib92,KT,Ksurvey,KL,Yam}. In addition,
the idea of \cite{BukhKlib} was extended to numerical methods for
coefficient inverse problems, see, e.g. \cite{Kconv}-\cite{KTnum}.

\textbf{Remark 1.1}. \emph{We are not concerned here with the issue of the
minimal smoothness. In doing so we follow the tradition of the field of
Inverse Problems, see, e.g. \cite{Nov1,Nov2}, \cite[Theorem 4.1]{Rom2}. }

We work below only with real valued functions. We formulate our two problems
in section 2. In section 3 we formulate Carleman estimates of \cite%
{MFG1,MFG2}. We prove H\"{o}lder stability estimates and uniqueness of our
two problems in sections 4 and 5.

\section{Two Problems}

\label{sec:2}

Below $\beta =const.>0.$ Let $x=\left( x_{1},x_{2},...,x_{n}\right) \in 
\mathbb{R}^{n}$\ denotes the position $x$\ of an agent and $t\geq 0$ denotes
time. Let $\Omega \subset \mathbb{R}^{n}$ be a bounded domain with the
piecewise smooth boundary $\partial \Omega $ and $T>0$ be a number. Denote 
\begin{equation*}
Q_{T}=\Omega \times \left( 0,T\right) ,S_{T}=\partial \Omega \times \left(
0,T\right) .
\end{equation*}

Recall that $u(x,t)$\ the value function and $m(x,t)$\ is the density of
players at the position $x$ and at the moment of time $t$. The conventional
MFGS of the second order consists of the system of two homogeneous nonlinear
parabolic PDEs with times running in two different directions \cite{A,LL3}.
However, we consider in this paper a more general case of two heterogeneous
parabolic PDEs. To do this, we incorporate non-zero terms in the right hand
sides to the conventional PDEs forming the MFGS. Hence, this is a
generalized MFGS of the form:%
\begin{equation}
\left. 
\begin{array}{c}
u_{t}(x,t)+\beta \Delta u(x,t){-r(x)(\nabla u(x,t))^{2}/2}+ \\ 
+F\left( \dint\limits_{\Omega }M\left( x,y\right) m\left( y,t\right)
dy,m\left( x,t\right) \right) =G_{1}\left( x,t\right) ,\text{ }\left(
x,t\right) \in Q_{T},%
\end{array}%
\right.  \label{2.1}
\end{equation}%
\begin{equation}
m_{t}(x,t)-\beta \Delta m(x,t){-\func{div}(r(x)m(x,t)\nabla u(x,t))}%
=G_{2}\left( x,t\right) ,\text{ }\left( x,t\right) \in Q_{T}.  \label{2.2}
\end{equation}

\bigskip Here, the coefficient $r(x)\in C^{1}\left( \overline{\Omega }%
\right) $ is similar with the elasticity of the medium, the function $F$ is
the interaction term. We assume zero Neumann boundary conditions, i.e. the
full reflection from the boundary%
\begin{equation}
\partial _{\nu }u\mid _{S_{T}}=\partial _{\nu }m\mid _{S_{T}}=0,  \label{2.3}
\end{equation}%
where $\nu =\nu \left( x\right) $ is the unit outward looking normal vector
at the point $\left( x,t\right) \in S_{T}.$ We consider in this paper the
following two problems:

\textbf{Problem 1}. \emph{Assuming that functions }$u,m\in H^{2}\left(
Q_{T}\right) $\emph{\ satisfy conditions (\ref{2.1})-(\ref{2.3}), obtain a H%
\"{o}lder stability estimate and uniqueness theorem for the case when the
following two functions }$u_{T}\left( x\right) $\emph{\ and }$m_{T}\left(
x\right) $\emph{\ are known:}%
\begin{equation}
u\left( x,T\right) =u_{T}\left( x\right) ,\text{ }m\left( x,T\right)
=m_{T}\left( x\right) ,\text{ }x\in \Omega .  \label{2.4}
\end{equation}

\textbf{Problem 2}. \emph{Assuming that functions }$u,m\in H^{2}\left(
Q_{T}\right) $\emph{\ satisfy conditions (\ref{2.1})-(\ref{2.3}), obtain a H%
\"{o}lder stability estimate and uniqueness theorem for the case when the
following two functions }$u_{0}\left( x\right) $\emph{\ and }$m_{0}\left(
x\right) $\emph{\ are known:}%
\begin{equation}
u\left( x,0\right) =u_{0}\left( x\right) ,\text{ }m\left( x,0\right)
=m_{0}\left( x\right) ,\text{ }x\in \Omega .  \label{2.5}
\end{equation}

\textbf{Remark 2.1}. \emph{The data in the right hand sides of (\ref{2.4})
and (\ref{2.5}) are measured with errors, so as the right hand sides of
equations (\ref{2.1}) and (\ref{2.2}). This is why stability estimates,
which are actually accuracy estimates for solutions of the MFGS, are
important in Problems 1,2, so as for problems considered in \cite%
{MFG1,MFG2,MFG3}.}

\section{Carleman Estimates}

\label{sec:3}

A Carleman estimate for a partial differential operator is always proven
only for the principal part of this operator since it is independent on its
lower order terms \cite[Lemma 2.1.1]{KL}. Therefore we formulate in this
section Carleman estimates for principal parts $\partial _{t}+\beta \Delta $%
, $\partial _{t}-\beta \Delta $ of operators of equations (\ref{2.1}), (\ref%
{2.2}). Carleman estimates for Problems 1 and 2 are different. The
difference is in the difference of Carleman Weight Functions, i.e. weight
functions involved in the resulting integral inequalities. Denote%
\begin{equation*}
H_{0}^{2}\left( Q_{T}\right) =\left\{ u\in H^{2}\left( Q_{T}\right)
:\partial _{\nu }u\mid _{S_{T}}=0\right\} .
\end{equation*}

\subsection{Carleman estimates for Problem 1}

\label{sec:3.1}

Introduce three parameters $b>0,\lambda >0$ and $k>2$. Also, introduce our
first Carleman Weight Function 
\begin{equation}
\varphi _{\lambda ,k}\left( t\right) =\exp \left( \lambda \left( t+b\right)
^{k}\right) ,t\in \left( 0,T\right) .  \label{3.1}
\end{equation}

\textbf{Theorem 3.1 }\cite{MFG1}\textbf{.} \emph{There exists a number }$%
C_{1}=C_{1}\left( b,T,\beta \right) >0,$ \emph{which depends only on listed
parameters, such that the following Carleman estimate is valid:}%
\begin{equation}
\left. 
\begin{array}{c}
\dint\limits_{Q_{T}}\left( u_{t}+\beta \Delta u\right) ^{2}\varphi _{\lambda
,k}^{2}dxdt\geq C_{1}\dint\limits_{Q_{T}}\left( u_{t}^{2}+\left( \Delta
u\right) ^{2}\right) \varphi _{\lambda ,k}^{2}dxdt+ \\ 
+C_{1}\lambda k\dint\limits_{Q_{T}}\left( \nabla u\right) ^{2}\varphi
_{\lambda ,k}^{2}dxdt+C_{1}\lambda ^{2}k^{2}\dint\limits_{Q_{T}}u^{2}\varphi
_{\lambda ,k}^{2}dxdt- \\ 
-e^{2\lambda \left( T+b\right) ^{k}}\dint\limits_{\Omega }\left[ \left(
\nabla _{x}u\right) ^{2}+\lambda k\left( T+b\right) ^{k}u^{2}\right] \left(
x,T\right) dx, \\ 
\forall \lambda >0,\forall k>2,\forall u\in H_{0}^{2}\left( Q_{T}\right) .%
\end{array}%
\right.  \label{3.2}
\end{equation}

\textbf{Theorem 3.2} \cite{MFG1}\textbf{.} \emph{There exist a sufficiently
large number }$k_{0}=k_{0}\left( \beta ,T,b\right) >2$\emph{\ and a number }$%
C=C\left( T,b\right) >0$\emph{\ depending only on listed parameters such
that the following Carleman estimate holds:} 
\begin{equation}
\left. 
\begin{array}{c}
\dint\limits_{Q_{T}}\left( u_{t}-\beta \Delta u\right) ^{2}\varphi _{\lambda
,k}^{2}dxdt\geq \\ 
\geq C_{1}\sqrt{k}\beta \dint\limits_{Q_{T}}\left( \nabla u\right)
^{2}\varphi _{\lambda ,k}^{2}dxdt+C_{1}\lambda
k^{2}\dint\limits_{Q_{T}}u^{2}\varphi _{\lambda ,k}^{2}dxdt- \\ 
-C_{1}\lambda k\left( T+b\right) ^{k-1}e^{2\lambda \left( T+b\right)
^{k}}\dint\limits_{\Omega }u^{2}\left( x,T\right) dx-C_{1}e^{2\lambda
b^{k}}\dint\limits_{\Omega }\left[ \left( \nabla u\right) ^{2}+\sqrt{\nu }%
u^{2}\right] \left( x,0\right) dx, \\ 
\forall \lambda >0,\forall k\geq k_{0}\left( \beta ,T,b\right) >2,\forall
u\in H_{0}^{2}\left( Q_{T}\right) ,%
\end{array}%
\right.  \label{3.3}
\end{equation}%
\emph{where the number }$C_{1}=C_{1}\left( a,T,\beta \right) >0$\emph{\
depends on the same parameters as ones in Theorem 3.1.}

In addition, we formulate a new integral identity, which was proven in \cite%
{MFG1}\textbf{. }In the past only a similar inequality rather than identity
was known \cite[Chapter 2, \S 6]{Lad}.

\textbf{Lemma 3.1.} \emph{Suppose that the domain }$\Omega $\emph{\ is a
rectangular prism. Then the following integral identity holds:}%
\begin{equation*}
\dint\limits_{\Omega }\left( \Delta u\right)
^{2}dx=\dsum\limits_{i,j=1}^{n}\dint\limits_{\Omega }u_{x_{i}x_{j}}^{2}dx,
\end{equation*}%
\begin{equation*}
\forall u\in \left\{ u\in H^{2}\left( \Omega \right) :\partial _{\nu }u\mid
_{\Omega }=0\right\} .
\end{equation*}

\subsection{ Carleman estimates for Problem 2}

\label{sec:3.2}

Let $c>2$ be a number.\ Let $\lambda >2$ be a sufficiently large parameter.
We will choose parameters $c$ and $\lambda $ later. Introduce the second
Carleman Weight Function $\varphi _{\lambda }\left( t\right) ,$%
\begin{equation}
\varphi _{\lambda }\left( t\right) =\exp \left( \left( T-t+c\right)
^{\lambda }\right) ,\text{ }t\in \left( 0,T\right) .  \label{3.4}
\end{equation}

\textbf{Theorem 3.3 }\cite{MFG2}.\ \emph{Choose the number }$c$\emph{\ in (%
\ref{3.4}) as: }$c>2.$\emph{\ Define the number }$\lambda _{0}$\emph{\ as: }%
\begin{equation}
\lambda _{0}=\lambda _{0}\left( T,c\right) =16\left( T+c\right)
^{2}>16c^{2}>64.  \label{3.7}
\end{equation}%
\emph{Then\ the following Carleman estimate is valid:}%
\begin{equation*}
\dint\limits_{Q_{T}}\left( u_{t}+\beta \Delta u\right) ^{2}\varphi _{\lambda
}^{2}dxdt\geq
\end{equation*}%
\begin{equation}
\geq C_{2}\sqrt{\lambda }\dint\limits_{Q_{T}}\left( \nabla u\right)
^{2}\varphi _{\lambda }^{2}dxdt+C_{2}\lambda ^{2}c^{\lambda
-2}\dint\limits_{Q_{T}}u^{2}\varphi _{\lambda }^{2}dxdt+  \label{3.8}
\end{equation}%
\begin{equation*}
-C_{2}e^{2c^{\lambda }}\dint\limits_{\Omega }\left( \left( \nabla u\right)
^{2}+u^{2}\right) \left( x,T\right) -C_{2}\lambda \left( T+c\right)
^{\lambda -1}e^{2\left( T+c\right) ^{\lambda }}\dint\limits_{\Omega
}u^{2}\left( x,0\right) dx,\text{ }
\end{equation*}%
\begin{equation*}
\forall \lambda \geq \lambda _{0},\forall u\in H_{0}^{2}\left( Q_{T}\right) ,
\end{equation*}%
\emph{where the constant }$C_{2}=C_{2}\left( c,T,\beta \right) >0$\emph{\
depends only on listed parameters.}

Theorem 3.4 is not exactly a Carleman estimate but rather a quasi-Carleman
estimate. This is because of two test functions $u$ and $v$ are involved in
it rather than just a single one.

\textbf{Theorem 3.4 }(a quasi-Carleman estimate) \cite{MFG2}. \emph{Let the
numbers }$c$\emph{\ and }$\lambda _{0}$\emph{\ be the same as the ones in
Theorem 3.3. Let the function }$g\in H^{1}\left( Q_{T}\right) $\emph{\ and }%
\begin{equation*}
\sup_{Q_{T}}\left\vert g\right\vert ,\sup_{Q_{T}}\left\vert \nabla
g\right\vert <\infty .
\end{equation*}%
\emph{\ Then the following quasi-Carleman estimate holds for any two
functions }$u,v\in H_{0}^{2}\left( Q_{T}\right) $\emph{:}%
\begin{equation*}
\dint\limits_{Q_{T}}\left( u_{t}-\beta \Delta u+g\Delta v\right) ^{2}\varphi
_{\lambda }^{2}\geq
\end{equation*}%
\begin{equation*}
\geq C_{3}\lambda c^{\lambda -1}\dint\limits_{Q_{T}}\left( \nabla u\right)
^{2}\varphi _{\lambda }^{2}dxdt+C_{3}\lambda ^{2}c^{2\lambda
-2}\dint\limits_{Q_{T}}u^{2}\varphi _{\lambda }^{2}dxdt-
\end{equation*}%
\begin{equation}
-C_{3}\lambda \left( T+c\right) ^{\lambda -1}\dint\limits_{Q_{T}}\left(
\nabla u\right) ^{2}\varphi _{\lambda }^{2}dxdt-  \label{3.9}
\end{equation}%
\begin{equation*}
-\lambda \left( T+c\right) ^{\lambda -1}e^{2\left( T+c\right) ^{\lambda
}}\dint\limits_{\Omega }u^{2}\left( x,0\right) dx,\text{ }\forall \lambda
\geq \lambda _{0},
\end{equation*}%
\emph{where the number }$C_{3}=C_{3}\left( \beta ,c,\left\Vert g\right\Vert
_{C^{1}\left( \overline{Q}_{T}\right) }\right) >0$\emph{\ depends only on
listed parameters.}

\section{H\"{o}lder Stability and Uniqueness for Problem 1}

\label{sec:4}

\bigskip For any number $\varepsilon \in \left( 0,T\right) $ define the
domain $Q_{\varepsilon ,T}$ as:%
\begin{equation}
Q_{\varepsilon ,T}=\Omega \times \left( \varepsilon ,T\right) \subset Q_{T}.
\label{4.0}
\end{equation}

\textbf{Theorem 4.1.} \emph{Let }$D_{1},D_{2},D_{3},D_{4}>0$\emph{\ be
certain numbers. Let in (\ref{2.1}), (\ref{2.2}) functions }$G_{1},G_{2}\in
L_{2}\left( Q_{T}\right) .$ \emph{Let }$F$\emph{\ be the function\ in (\ref%
{2.1}). Assume that the function }$F\left( y,z\right) :\mathbb{R}%
^{2}\rightarrow \mathbb{R}$\emph{\ has derivatives }$F_{y},F_{z}\in C\left( 
\mathbb{R}^{2}\right) $ \emph{such that} 
\begin{equation}
\max \left( \sup_{\mathbb{R}^{2}}\left\vert F_{y}\left( y,z\right)
\right\vert ,\sup_{\mathbb{R}^{2}}\left\vert F_{y}\left( y,z\right)
\right\vert \right) \leq D_{1}.  \label{4.1}
\end{equation}%
\emph{In (\ref{2.1}), (\ref{2.2}), let the functions }$M\left( x,y\right) $%
\emph{\ and }$r\in C^{1}\left( \overline{\Omega }\right) $ be such that 
\begin{equation}
\sup_{\Omega \times \Omega }\left\vert M\left( x,y\right) \right\vert
,\left\Vert r\right\Vert _{C^{1}\left( \overline{\Omega }\right) }\leq D_{2}.
\label{4.2}
\end{equation}%
\emph{Define sets of functions }$K_{3}\left( D_{3}\right) $\emph{\ and }$%
K_{4}\left( D_{4}\right) $\emph{\ as}%
\begin{equation}
K_{3}\left( D_{3}\right) =\left\{ u\in H_{0}^{2}\left( Q_{T}\right)
:\sup_{Q_{T}}\left\vert u\right\vert ,\sup_{Q_{T}}\left\vert \nabla
u\right\vert ,\sup_{Q_{T}}\left\vert \Delta u\right\vert \leq D_{3}\right\} ,
\label{4.3}
\end{equation}%
\begin{equation}
K_{4}\left( D_{4}\right) =\left\{ u\in H_{0}^{2}\left( Q_{T}\right)
:\sup_{Q_{T}}\left\vert u\right\vert ,\sup_{Q_{T}}\left\vert \nabla
u\right\vert \leq D_{4}\right\} .  \label{4.4}
\end{equation}%
\emph{Let}%
\begin{equation}
D=\max \left( D_{1},D_{2},D_{3},D_{4}\right) .  \label{4.5}
\end{equation}%
\emph{Assume that two pairs of functions }$\left( u_{1},m_{1}\right) $\emph{%
\ and }$\left( u_{2},m_{2}\right) $\emph{\ satisfy equations (\ref{2.1})- (%
\ref{2.3}) with two pairs of functions }$\left( G_{1,1},G_{2,1}\right) $%
\emph{\ and }$\left( G_{1,2},G_{2,2}\right) $\emph{\ respectively and are
such that }%
\begin{equation}
\left( u_{1},m_{1}\right) ,\left( u_{2},m_{2}\right) \in K_{3}\left(
D_{3}\right) \times K_{4}\left( D_{4}\right) .  \label{4.6}
\end{equation}%
\emph{\ Assume that these two pairs of functions }$\left( u_{1},m_{1}\right) 
$\emph{\ and }$\left( u_{2},m_{2}\right) $ \emph{have the following terminal
conditions (see (\ref{2.4})):}%
\begin{equation}
u_{1}\left( x,T\right) =u_{T}^{\left( 1\right) }\left( x\right) ,\text{ }%
m_{1}\left( x,T\right) =m_{T}^{\left( 1\right) }\left( x\right) ,\text{ }%
x\in \Omega ,  \label{4.7}
\end{equation}%
\begin{equation}
\text{ }u_{2}\left( x,T\right) =u_{T}^{\left( 2\right) }\left( x\right) ,%
\text{ }m_{2}\left( x,T\right) =m_{T}^{\left( 2\right) }\left( x\right) ,%
\text{ }x\in \Omega .  \label{4.8}
\end{equation}%
\emph{Let the number }$\varepsilon \in \left( 0,T\right) $\emph{\ and let }$%
Q_{\varepsilon ,T}$\emph{\ be the domain defined in (\ref{4.0}). Then there
exists a number }$C_{4}=C_{4}\left( \beta ,D,T,\Omega ,\varepsilon \right)
>0 $\emph{\ and a sufficiently small number }$\delta _{0}=$\emph{\ }$\delta
_{0}\left( \beta ,D,T,\Omega ,\varepsilon \right) \in \left( 0,1\right) $ 
\emph{depending only on listed parameters such that if \ }$\delta \in \left(
0,\delta _{0}\right) $ \emph{and} \emph{\ }%
\begin{equation}
\left\Vert u_{T}^{\left( 1\right) }-u_{T}^{\left( 2\right) }\right\Vert
_{H^{1}\left( \Omega \right) },\left\Vert m_{T}^{\left( 1\right)
}-m_{T}^{\left( 2\right) }\right\Vert _{L_{2}\left( \Omega \right) }\leq
\delta ,  \label{4.9}
\end{equation}%
\begin{equation}
\left\Vert G_{1,1}-G_{1,2}\right\Vert _{L_{2}\left( Q_{T}\right)
},\left\Vert G_{2,1}-G_{2,2}\right\Vert _{L_{2}\left( Q_{T}\right) }\leq
\delta ,  \label{4.90}
\end{equation}%
\emph{then there exists a number }$\rho =\rho \left( \beta ,D,T,\Omega
,\varepsilon \right) \in \left( 0,1/6\right) $\emph{\ depending only on
listed parameters such that the following two H\"{o}lder stability estimates
are valid:}%
\begin{equation}
\left. 
\begin{array}{c}
\left\Vert \partial _{t}u_{1}-\partial _{t}u_{2}\right\Vert _{L_{2}\left(
Q_{\varepsilon ,T}\right) }+\left\Vert \Delta u_{1}-\Delta u_{2}\right\Vert
_{L_{2}\left( Q_{\varepsilon ,T}\right) }+\left\Vert u_{1}-u_{2}\right\Vert
_{H^{1,0}\left( Q_{\varepsilon ,T}\right) }\leq \\ 
\leq C_{4}\left( 1+\left\Vert m_{1}-m_{2}\right\Vert _{H^{2}\left(
Q_{T}\right) }\right) \delta ^{\rho },\text{ }\forall \delta \in \left(
0,\delta _{0}\right) ,%
\end{array}%
\right.  \label{4.10}
\end{equation}%
\begin{equation}
\left\Vert m_{1}-m_{2}\right\Vert _{H^{1,0}\left( Q_{_{\varepsilon
,}T}\right) }\leq C_{4}\left( 1+\left\Vert m_{1}-m_{2}\right\Vert
_{H^{2}\left( Q_{T}\right) }\right) \delta ^{\rho }.  \label{4.11}
\end{equation}%
\emph{\ Furthermore, if the domain }$\Omega $\emph{\ is a rectangular prism,
then estimate (\ref{4.10}) can be strengthened as:} 
\begin{equation}
\left\Vert u_{1}-u_{2}\right\Vert _{H^{2,1}\left( Q_{\varepsilon ,T}\right)
}\leq C_{4}\left( 1+\left\Vert m_{1}-m_{2}\right\Vert _{H^{2}\left(
Q_{T}\right) }\right) \delta ^{\rho }.  \label{4.12}
\end{equation}%
\emph{Furthermore, if in (\ref{4.7}) and (\ref{4.8}) }%
\begin{equation}
u_{T}^{\left( 1\right) }\left( x\right) \equiv u_{T}^{\left( 2\right)
}\left( x\right) ,\text{ }m_{T}^{\left( 1\right) }\left( x\right) \equiv
m_{T}^{\left( 2\right) }\left( x\right) ,\text{ }x\in \Omega ,  \label{4.13}
\end{equation}%
\begin{equation}
G_{1,1}\left( x,t\right) \equiv G_{1,2}\left( x,t\right) ,\text{ }%
G_{2,1}\left( x,t\right) \equiv G_{2,2}\left( x,t\right) ,\text{ }\left(
x,t\right) \in Q_{T},\text{ }  \label{4.130}
\end{equation}%
\emph{then }$u_{1}\left( x,t\right) \equiv u_{2}\left( x,t\right) $\emph{\
and }$m_{1}\left( x,t\right) \equiv m_{2}\left( x,t\right) $ \emph{in }$%
Q_{T},$\emph{\ which means that Problem 1 has at most one solution }$\left(
u,m\right) \in K_{3}\left( D_{3}\right) \times K_{4}\left( D_{4}\right) $%
\emph{.}

\textbf{Remark 4.1: }\emph{Consider the condition of this theorem that two
pairs }$\left( u_{1},m_{1}\right) ,\left( u_{2},m_{2}\right) $\emph{\ belong
to an a priori chosen bounded set }$K_{3}\left( D_{3}\right) \times
K_{4}\left( D_{4}\right) .$\emph{\ Such conditions are typical ones in the
theory of ill-posed problems, see, e.g. \cite{BK,LRS}. }

\textbf{Proof of Theorem 4.1.} In this proof $\widetilde{C}_{4}=\widetilde{C}%
_{4}\left( \beta ,D,T,\Omega \right) >0$ denotes different numbers depending
only on $\beta ,D,T,\Omega $ and $C_{4}=C_{4}\left( \beta ,D,T,\Omega
,\varepsilon \right) >0$ denotes different numbers depending not only on
parameters $\beta ,D,T,\Omega $ but on $\varepsilon $ as well. Consider four
arbitrary numbers $y_{1},z_{1},y_{2},z_{2}\in \mathbb{R}.$ Let $\widetilde{y}%
=y_{1}-y_{2}$ and $\widetilde{z}=z_{1}-z_{2}.$ Hence, 
\begin{equation}
y_{1}z_{1}-y_{2}z_{2}=\widetilde{y}z_{1}+\widetilde{z}y_{2}.  \label{4.14}
\end{equation}%
Denote%
\begin{equation}
v\left( x,t\right) =u_{1}\left( x,t\right) -u_{2}\left( x,t\right) ,\text{ }%
p\left( x,t\right) =m_{1}\left( x,t\right) -m_{2}\left( x,t\right) ,\left(
x,t\right) \in Q_{T},  \label{4.15}
\end{equation}%
\begin{equation}
v_{T}\left( x\right) =u_{T}^{\left( 1\right) }\left( x\right) -u_{T}^{\left(
2\right) }\left( x\right) ,\text{ }p_{T}\left( x\right) =m_{T}^{\left(
1\right) }\left( x\right) -m_{T}^{\left( 2\right) }\left( x\right) ,\text{ }%
x\in \Omega ,  \label{4.16}
\end{equation}%
\begin{equation}
\widetilde{G}_{1}\left( x,t\right) =\left( G_{1,1}-G_{1,2}\right) \left(
x,t\right) ,\text{ }\widetilde{G}_{2}\left( x,t\right) =\left(
G_{2,1}-G_{2,2}\right) \left( x,t\right) ,\left( x,t\right) \in Q_{T}.
\label{4.160}
\end{equation}

Using (\ref{4.1})-(\ref{4.6}) and the multidimensional analog of Taylor
formula \cite{V}, we obtain%
\begin{equation}
\left. 
\begin{array}{c}
F\left( \dint\limits_{\Omega }M\left( x,y\right) m_{1}\left( y,t\right)
dy,m_{1}\left( x,t\right) \right) - \\ 
-F\left( \dint\limits_{\Omega }M\left( x,y\right) m_{2}\left( y,t\right)
dy,m_{2}\left( x,t\right) \right) = \\ 
=f_{1}\left( x,t\right) \dint\limits_{\Omega }M\left( x,y\right) p\left(
y,t\right) +f_{2}\left( x,t\right) p\left( x,t\right) ,%
\end{array}%
\right.  \label{4.17}
\end{equation}%
where functions $f_{1},f_{2}$ are such that 
\begin{equation}
\left\vert f_{1}\left( x,t\right) \right\vert ,\left\vert f_{2}\left(
x,t\right) \right\vert \leq D,\text{ }\left( x,t\right) \in Q_{T}.
\label{4.18}
\end{equation}

Subtract equations (\ref{2.1}), (\ref{2.2}) for the pair $\left(
u_{2},m_{2}\right) $ from corresponding equations for the pair $\left(
u_{1},m_{1}\right) $. Then use (\ref{4.6})-(\ref{4.8}), (\ref{4.14})-(\ref%
{4.18}) and recall that Carleman estimates can work with both equations and
inequalities \cite{BK,KT,KL,LRS,Yam}. Hence, it is convenient to replace
resulting equations with two inequalities:%
\begin{equation}
\left\vert v_{t}+\beta \Delta v\right\vert \left( x,t\right) \leq \widetilde{%
C}_{4}\left( \left\vert \nabla v\right\vert +\dint\limits_{\Omega
}\left\vert p\left( y,t\right) \right\vert dy+\left\vert p\right\vert
+\left\vert \widetilde{G}_{1}\right\vert \right) \left( x,t\right) ,\text{ }%
\left( x,t\right) \in Q_{T},  \label{4.19}
\end{equation}%
\begin{equation}
\left\vert p_{t}-\beta \Delta p\right\vert \left( x,t\right) \leq \widetilde{%
C}_{4}\left( \left\vert \nabla p\right\vert +\left\vert p\right\vert
+\left\vert \Delta v\right\vert +\left\vert \nabla v\right\vert +\left\vert 
\widetilde{G}_{2}\right\vert \right) \left( x,t\right) ,\text{ }\left(
x,t\right) \in Q_{T},  \label{4.20}
\end{equation}%
\begin{equation}
\partial _{\nu }v\mid _{S_{T}}=0,\text{ }\partial _{\nu }p\mid _{S_{T}}=0,
\label{4.21}
\end{equation}%
\begin{equation}
v\left( x,T\right) =v_{T}\left( x\right) ,\text{ }p\left( x,T\right)
=p_{T}\left( x\right) ,\text{ }x\in \Omega .  \label{4.22}
\end{equation}

Squaring both sides of equation (\ref{4.19}) and (\ref{4.20}), applying
Cauchy-Schwarz inequality, multiplying by the CWF $\varphi _{\lambda
,k}^{2}\left( t\right) $ defined in (\ref{3.1}) and integrating over $Q_{T},$
we obtain%
\begin{equation}
\left. 
\begin{array}{c}
\dint\limits_{Q_{T}}\left( v_{t}+\beta \Delta v\right) ^{2}\varphi _{\lambda
,k}^{2}dxdt\leq \widetilde{C}_{4}\dint\limits_{Q_{T}}\left( \nabla v\right)
^{2}\varphi _{\lambda ,k}^{2}dxdt+ \\ 
+\widetilde{C}_{4}\dint\limits_{Q_{T}}\left( p^{2}+\dint\limits_{\Omega
}p^{2}\left( y,t\right) dy+\widetilde{G}_{1}^{2}\right) \varphi _{\lambda
,k}^{2}dxdt,%
\end{array}%
\right.  \label{4.23}
\end{equation}%
\begin{equation}
\left. 
\begin{array}{c}
\dint\limits_{Q_{T}}\left( p_{t}-\beta \Delta p\right) ^{2}\varphi _{\lambda
,k}^{2}dxdt\leq \widetilde{C}_{4}\dint\limits_{Q_{T}}\left( \left( \nabla
p\right) ^{2}+p^{2}\right) \varphi _{\lambda ,k}^{2}dxdt+ \\ 
+\widetilde{C}_{4}\dint\limits_{Q_{T}}\left( \left( \Delta v\right)
^{2}+\left( \nabla v\right) ^{2}+v^{2}+\widetilde{G}_{2}^{2}\right) \varphi
_{\lambda ,k}^{2}dxdt.%
\end{array}%
\right.  \label{4.24}
\end{equation}%
Note that 
\begin{equation}
\dint\limits_{Q_{T}}\left( \dint\limits_{\Omega }p^{2}\left( y,t\right)
dy\right) \varphi _{\lambda ,k}^{2}\left( t\right) dxdt\leq \widetilde{C}%
_{4}\dint\limits_{Q_{T}}p^{2}\varphi _{\lambda ,k}^{2}\left( t\right) dxdt.
\label{4.25}
\end{equation}

Set 
\begin{equation}
b=1  \label{4.250}
\end{equation}%
in the Carleman Weight Function $\varphi _{\lambda ,k}\left( t\right) .$
Apply Carleman estimate (\ref{3.2}) to the left hand side of (\ref{4.23})
and use (\ref{4.21}), (\ref{4.22}) and (\ref{4.25}). We obtain 
\begin{equation}
\left. 
\begin{array}{c}
\widetilde{C}_{4}\dint\limits_{Q_{T}}\left( p^{2}+\widetilde{G}%
_{1}^{2}\right) \varphi _{\lambda ,k}^{2}dxdt\geq \\ 
\geq \dint\limits_{Q_{T}}\left( v_{t}^{2}+\left( \Delta v\right) ^{2}\right)
\varphi _{\lambda ,k}^{2}dxdt+ \\ 
+\lambda k\dint\limits_{Q_{T}}\left( \nabla v\right) ^{2}\varphi _{\lambda
,k}^{2}dxdt+\lambda ^{2}k^{2}\dint\limits_{Q_{T}}v^{2}\varphi _{\lambda
,k}^{2}dxdt- \\ 
-\widetilde{C}_{4}e^{2\lambda \left( T+1\right) ^{k}}\dint\limits_{\Omega }%
\left[ \left( \nabla _{x}v_{T}\right) ^{2}+\lambda k\left( T+1\right)
^{k}v_{T}^{2}\right] dx, \\ 
\forall \lambda >0,\forall k>2.%
\end{array}%
\right.  \label{4.26}
\end{equation}%
Choosing $\lambda _{1}=\lambda _{1}\left( \beta ,D,T,\Omega \right) \geq 1$
so large that 
\begin{equation}
\lambda _{1}>2\widetilde{C}_{4}  \label{4.270}
\end{equation}%
and recalling that $k>2,$ we obtain from (\ref{4.26}) 
\begin{equation}
\left. 
\begin{array}{c}
\widetilde{C}_{4}\dint\limits_{Q_{T}}p^{2}\varphi _{\lambda ,k}^{2}dxdt+%
\widetilde{C}_{4}\dint\limits_{Q_{T}}\widetilde{G}_{1}^{2}\varphi _{\lambda
,k}^{2}dxdt\geq \\ 
\geq \dint\limits_{Q_{T}}\left( v_{t}^{2}+\left( \Delta v\right) ^{2}\right)
\varphi _{\lambda ,k}^{2}dxdt+ \\ 
+\lambda k\dint\limits_{Q_{T}}\left( \nabla v\right) ^{2}\varphi _{\lambda
,k}^{2}dxdt+\lambda ^{2}k^{2}\dint\limits_{Q_{T}}v^{2}\varphi _{\lambda
,k}^{2}dxdt- \\ 
-\widetilde{C}_{4}e^{2\lambda \left( T+1\right) ^{k}}\dint\limits_{\Omega }%
\left[ \left( \nabla _{x}v_{T}\right) ^{2}+\lambda k\left( T+1\right)
^{k}v_{T}^{2}\right] dx, \\ 
\forall \lambda \geq \lambda _{1},\forall k>2.%
\end{array}%
\right.  \label{4.28}
\end{equation}

We now apply Carleman estimate (\ref{3.3}) to the left hand side of (\ref%
{4.24}). We obtain%
\begin{equation}
\left. 
\begin{array}{c}
\widetilde{C}_{4}\dint\limits_{Q_{T}}\left( \left( \Delta v\right)
^{2}+v^{2}+\widetilde{G}_{2}^{2}\right) \varphi _{\lambda ,k}^{2}dxdt+%
\widetilde{C}_{4}\dint\limits_{Q_{T}}\left( \left( \nabla p\right)
^{2}+p^{2}\right) \varphi _{\lambda ,k}^{2}dxdt\geq \\ 
\geq \sqrt{k}\dint\limits_{Q_{T}}\left( \nabla p\right) ^{2}\varphi
_{\lambda ,k}^{2}dxdt+\lambda k^{2}\dint\limits_{Q_{T}}p^{2}\varphi
_{\lambda ,k}^{2}dxdt- \\ 
-\lambda k\left( T+b\right) ^{k-1}e^{2\lambda \left( T+b\right)
^{k}}\dint\limits_{\Omega }p_{T}^{2}\left( x\right) dx-e^{2\lambda
b^{k}}\dint\limits_{\Omega }\left[ \left( \nabla p\right) ^{2}+\sqrt{k}p^{2}%
\right] \left( x,0\right) dx, \\ 
\forall \lambda >0,\forall k\geq k_{0}=k_{0}\left( \beta ,T\right) >2.%
\end{array}%
\right.  \label{4.29}
\end{equation}%
Choose the number $k_{0}=k_{0}\left( \beta ,T\right) $ so large that 
\begin{equation}
k_{0}>\max \left( 2,4\widetilde{C}_{4}^{2}\right) ,  \label{4.290}
\end{equation}%
and, until (\ref{4.410}), set $k=k_{0}.$ \ Also, let $\lambda \geq \lambda
_{1},$ where the number $\lambda _{1}=\lambda _{1}\left( \beta ,D,T,\Omega
\right) \geq 1$ is defined in (\ref{4.27}). Then (\ref{4.29}) implies%
\begin{equation}
\left. 
\begin{array}{c}
\widetilde{C}_{4}\dint\limits_{Q_{T}}\left( \left( \Delta v\right)
^{2}+v^{2}+\widetilde{G}_{2}^{2}\right) \varphi _{\lambda ,k}^{2}dxdt\geq \\ 
\geq \dint\limits_{Q_{T}}\left( \nabla p\right) ^{2}\varphi _{\lambda
,k}^{2}dxdt+\lambda \dint\limits_{Q_{T}}p^{2}\varphi _{\lambda ,k}^{2}dxdt-
\\ 
-\widetilde{C}_{4}\lambda \left( T+1\right) ^{k-1}e^{2\lambda \left(
T+1\right) ^{k}}\dint\limits_{\Omega }p_{T}^{2}\left( x\right) dx-\widetilde{%
C}_{4}e^{2\lambda }\dint\limits_{\Omega }\left[ \left( \nabla p\right)
^{2}+p^{2}\right] \left( x,0\right) dx, \\ 
\forall \lambda \geq \lambda _{1}.%
\end{array}%
\right.  \label{4.30}
\end{equation}%
In particular, it follows from (\ref{4.30}) that%
\begin{equation}
\left. 
\begin{array}{c}
\dint\limits_{Q_{T}}p^{2}\varphi _{1,\lambda ,k}^{2}dxdt\leq \widetilde{C}%
_{4}\lambda ^{-1}\dint\limits_{Q_{T}}\left( \left( \Delta v\right)
^{2}+v^{2}+\widetilde{G}_{2}^{2}\right) \varphi _{\lambda ,k}^{2}dxdt+ \\ 
+\widetilde{C}_{4}e^{2\lambda \left( T+1\right) ^{k}}\dint\limits_{\Omega
}p_{T}^{2}\left( x\right) dx+\widetilde{C}_{4}e^{2\lambda
}\dint\limits_{\Omega }\left[ \left( \nabla p\right) ^{2}+p^{2}\right]
\left( x,0\right) dx,\text{ }\forall \lambda \geq \lambda _{1}.%
\end{array}%
\right.  \label{4.31}
\end{equation}%
Replacing the first term of the first line of (\ref{4.28}) with the right
hand side of inequality (\ref{4.31}), we obtain%
\begin{equation*}
\widetilde{C}_{4}\lambda e^{2\lambda \left( T+1\right)
^{k}}\dint\limits_{\Omega }\left[ p_{T}^{2}+\left( \nabla _{x}v_{T}\right)
^{2}+v_{T}^{2}\right] \left( x\right) dx+\widetilde{C}_{4}\dint%
\limits_{Q_{T}}\left( \widetilde{G}_{1}^{2}+\widetilde{G}_{2}^{2}\right)
\varphi _{\lambda ,k}^{2}dxdt+
\end{equation*}%
\begin{equation*}
+\widetilde{C}_{4}e^{2\lambda }\dint\limits_{\Omega }\left[ \left( \nabla
p\right) ^{2}+p^{2}\right] \left( x,0\right) dx+
\end{equation*}%
\begin{equation}
+\frac{\widetilde{C}_{4}}{\lambda }\dint\limits_{Q_{T}}\left( \left( \Delta
v\right) ^{2}+v^{2}\right) \varphi _{\lambda ,k}^{2}dxdt\geq  \label{4.32}
\end{equation}%
\begin{equation*}
\geq \dint\limits_{Q_{T}}\left( v_{t}^{2}+\left( \Delta v\right) ^{2}\right)
\varphi _{\lambda ,k}^{2}dxdt+
\end{equation*}%
\begin{equation*}
+\lambda \dint\limits_{Q_{T}}\left( \nabla v\right) ^{2}\varphi _{\lambda
,k}^{2}dxdt+\lambda ^{2}\dint\limits_{Q_{T}}v^{2}\varphi _{\lambda
,k}^{2}dxdt,\text{ }\forall \lambda \geq \lambda _{1}.
\end{equation*}%
By (\ref{4.270}) $\widetilde{C}_{4}/\lambda <1/2,$ $\forall \lambda \geq
\lambda _{1}.$ Hence, terms in the 4$^{\text{th}}$ and 5$^{\text{th}}$ lines
of (\ref{4.32}) absorb terms in the 3$^{\text{rd}}$ line of (\ref{4.32}).
Hence,%
\begin{equation*}
\widetilde{C}_{4}\lambda e^{2\lambda \left( T+1\right)
^{k}}\dint\limits_{\Omega }\left[ p_{T}^{2}+\left( \nabla _{x}v_{T}\right)
^{2}+v_{T}^{2}\right] \left( x\right) dx+\widetilde{C}_{4}\dint%
\limits_{Q_{T}}\left( \widetilde{G}_{1}^{2}+\widetilde{G}_{2}^{2}\right)
\varphi _{\lambda ,k}^{2}dxdt+
\end{equation*}%
\begin{equation*}
+\widetilde{C}_{4}e^{2\lambda }\dint\limits_{\Omega }\left[ \left( \nabla
p\right) ^{2}+p^{2}\right] \left( x,0\right) dx\geq
\end{equation*}%
\begin{equation}
\geq \dint\limits_{Q_{T}}\left( v_{t}^{2}+\left( \Delta v\right) ^{2}\right)
\varphi _{\lambda ,k}^{2}dxdt+  \label{4.33}
\end{equation}%
\begin{equation*}
+\lambda \dint\limits_{Q_{T}}\left( \nabla v\right) ^{2}\varphi _{\lambda
,k}^{2}dxdt+\lambda ^{2}\dint\limits_{Q_{T}}v^{2}\varphi _{\lambda
,k}^{2}dxdt,\text{ }\forall \lambda \geq \lambda _{1}.
\end{equation*}

Comparing the last two lines of (\ref{4.33}) with the first line of (\ref%
{4.30}), we obtain%
\begin{equation*}
\widetilde{C}_{4}\lambda e^{2\lambda \left( T+1\right)
^{k}}\dint\limits_{\Omega }\left[ p_{T}^{2}+\left( \nabla _{x}v_{T}\right)
^{2}+v_{T}^{2}\right] \left( x\right) dx+\widetilde{C}_{4}\dint%
\limits_{Q_{T}}\left( \widetilde{G}_{1}^{2}+\widetilde{G}_{2}^{2}\right)
\varphi _{\lambda ,k}^{2}dxdt+
\end{equation*}%
\begin{equation}
+\widetilde{C}_{4}e^{2\lambda }\dint\limits_{\Omega }\left[ \left( \nabla
p\right) ^{2}+p^{2}\right] \left( x,0\right) dx\geq  \label{4.34}
\end{equation}%
\begin{equation*}
\geq \dint\limits_{Q_{T}}\left( \nabla p\right) ^{2}\varphi _{\lambda
,k}^{2}dxdt+\lambda \dint\limits_{Q_{T}}p^{2}\varphi _{\lambda ,k}^{2}dxdt,%
\text{ }\forall \lambda \geq \lambda _{1}.
\end{equation*}%
Summing up (\ref{4.33}) and (\ref{4.34}), we obtain%
\begin{equation*}
\dint\limits_{Q_{T}}\left( v_{t}^{2}+\left( \Delta v\right) ^{2}+\left(
\nabla v\right) ^{2}+v^{2}+\left( \nabla p\right) ^{2}+p^{2}\right) \varphi
_{\lambda ,k}^{2}dxdt\leq
\end{equation*}%
\begin{equation}
\leq \widetilde{C}_{4}e^{3\lambda \left( T+1\right)
^{k}}\dint\limits_{\Omega }\left[ p_{T}^{2}+\left( \nabla _{x}v_{T}\right)
^{2}+v_{T}^{2}\right] \left( x\right) dx  \label{4.35}
\end{equation}%
\begin{equation*}
+\widetilde{C}_{4}e^{2\lambda }\dint\limits_{\Omega }\left[ \left( \nabla
p\right) ^{2}+p^{2}\right] \left( x,0\right) dx+\widetilde{C}%
_{4}\dint\limits_{Q_{T}}\left( \widetilde{G}_{1}^{2}+\widetilde{G}%
_{2}^{2}\right) \varphi _{\lambda ,k}^{2}dxdt,\text{ }\forall \lambda \geq
\lambda _{1}.
\end{equation*}%
Since by (\ref{4.0}) $Q_{\varepsilon ,T}\subset Q_{T},$ then replacing $%
Q_{T} $ with $Q_{\varepsilon ,T}$ in the first line of (\ref{4.35}), we
strengthen this inequality. Hence, 
\begin{equation*}
\dint\limits_{Q_{\varepsilon ,T}}\left( v_{t}^{2}+\left( \Delta v\right)
^{2}+\left( \nabla v\right) ^{2}+v^{2}+\left( \nabla p\right)
^{2}+p^{2}\right) \varphi _{1,\lambda ,k}^{2}dxdt\leq
\end{equation*}%
\begin{equation}
\leq \widetilde{C}_{4}e^{3\lambda \left( T+1\right)
^{k}}\dint\limits_{\Omega }\left[ p_{T}^{2}+\left( \nabla _{x}v_{T}\right)
^{2}+v_{T}^{2}\right] \left( x\right) dx  \label{4.36}
\end{equation}%
\begin{equation*}
+\widetilde{C}_{4}e^{2\lambda }\dint\limits_{\Omega }\left[ \left( \nabla
p\right) ^{2}+p^{2}\right] \left( x,0\right) dx+\widetilde{C}%
_{4}\dint\limits_{Q_{T}}\left( \widetilde{G}_{1}^{2}+\widetilde{G}%
_{2}^{2}\right) \varphi _{\lambda ,k}^{2}dxdt,\text{ }\forall \lambda \geq
\lambda _{1}.
\end{equation*}

Next, by (\ref{3.1}), (\ref{4.0}) and (\ref{4.250})%
\begin{equation}
\min_{\overline{Q}_{\varepsilon ,T}}\varphi _{\lambda ,k}^{2}\left( t\right)
=e^{2\lambda \left( \varepsilon +1\right) ^{k}},\text{ }  \label{4.37}
\end{equation}%
\begin{equation}
\max_{\overline{Q}_{T}}\varphi _{\lambda ,k}^{2}\left( t\right) =e^{2\lambda
\left( T+1\right) ^{k}}.  \label{4.38}
\end{equation}%
Also, by (\ref{4.9}), (\ref{4.90}), (\ref{4.16}), (\ref{4.160}) and (\ref%
{4.38}) 
\begin{equation*}
\widetilde{C}_{4}e^{3\lambda \left( T+1\right) ^{k}}\dint\limits_{\Omega }%
\left[ p_{T}^{2}+\left( \nabla _{x}v_{T}\right) ^{2}+v_{T}^{2}\right] \left(
x\right) dx+\widetilde{C}_{4}\dint\limits_{Q_{T}}\left( \widetilde{G}%
_{1}^{2}+\widetilde{G}_{2}^{2}\right) \varphi _{\lambda ,k}^{2}dxdt\leq
\end{equation*}%
\begin{equation}
\leq \widetilde{C}_{4}e^{3\lambda \left( T+1\right) ^{k}}\delta ^{2}.
\label{4.39}
\end{equation}%
By the trace theorem%
\begin{equation}
\widetilde{C}_{4}e^{2\lambda }\dint\limits_{\Omega }\left[ \left( \nabla
p\right) ^{2}+p^{2}\right] \left( x,0\right) dx\leq \widetilde{C}%
_{4}e^{2\lambda }\left\Vert p\right\Vert _{H^{2}\left( Q_{T}\right) }^{2}.
\label{4.40}
\end{equation}%
Hence, using (\ref{4.36}), (\ref{4.37}), (\ref{4.39}) and (\ref{4.40}), we
replace $\widetilde{C}_{4}$ with $C_{4}$ (see the beginning of this proof)
and obtain for all $\lambda \geq \lambda _{1}:$%
\begin{equation}
\left\Vert v_{t}\right\Vert _{L_{2}\left( Q_{\varepsilon ,T}\right)
}^{2}+\left\Vert \Delta v\right\Vert _{L_{2}\left( Q_{\varepsilon ,T}\right)
}^{2}+\left\Vert v\right\Vert _{H^{1,0}\left( Q_{\varepsilon ,T}\right)
}^{2}+\left\Vert p\right\Vert _{H^{1,0}\left( Q_{\varepsilon ,T}\right)
}^{2}\leq  \label{4.41}
\end{equation}%
\begin{equation*}
\leq C_{4}e^{3\lambda \left( T+1\right) ^{k}}\delta ^{2}+C_{4}\exp \left[
-2\lambda \left( \varepsilon +1\right) ^{k}\left( 1-\frac{1}{\left(
\varepsilon +1\right) ^{k}}\right) \right] \left\Vert p\right\Vert
_{H^{2}\left( Q_{T}\right) }^{2},
\end{equation*}%
Recalling (\ref{4.290}), choose $k_{1}=k_{1}\left( \beta ,T,\varepsilon
\right) \geq k_{0}\left( \beta ,T\right) $ so large that 
\begin{equation}
\frac{1}{\left( \varepsilon +1\right) ^{k_{1}}}<\frac{1}{2}  \label{4.410}
\end{equation}%
and set $k=k_{1}.$ Then (\ref{4.41}) implies%
\begin{equation}
\left\Vert v_{t}\right\Vert _{L_{2}\left( Q_{\varepsilon ,T}\right)
}^{2}+\left\Vert \Delta v\right\Vert _{L_{2}\left( Q_{\varepsilon ,T}\right)
}^{2}+\left\Vert v\right\Vert _{H^{1,0}\left( Q_{\varepsilon ,T}\right)
}^{2}+\left\Vert p\right\Vert _{H^{1,0}\left( Q_{\varepsilon ,T}\right)
}^{2}\leq  \label{4.42}
\end{equation}%
\begin{equation*}
\leq C_{4}e^{3\lambda \left( T+1\right) ^{k}}\delta ^{2}+C_{4}e^{-\lambda
\left( \varepsilon +1\right) ^{k}}\left\Vert p\right\Vert _{H^{2}\left(
Q_{T}\right) }^{2},\text{ }\lambda \geq \lambda _{1}.
\end{equation*}%
Choose $\lambda =\lambda \left( \delta \right) $ such that 
\begin{equation}
e^{3\lambda \left( \delta \right) \left( T+1\right) ^{k}}\delta ^{2}=\delta .
\label{4.43}
\end{equation}%
Hence, 
\begin{equation}
\lambda \left( \delta \right) =\ln \left[ \delta ^{\left( 3\left( T+1\right)
\right) ^{-1}}\right] ,  \label{4.44}
\end{equation}%
\begin{equation}
e^{-\lambda \left( \varepsilon +1\right) ^{k}}=\delta ^{2\rho },\text{ }%
2\rho =\frac{1}{3}\left( \frac{\varepsilon +1}{T+1}\right) ^{k}<\frac{1}{3}.
\label{4.45}
\end{equation}%
Choose $\delta _{0}=\delta _{0}\left( \beta ,D,T,\Omega ,\varepsilon \right)
\in \left( 0,1\right) $ so small that%
\begin{equation}
\lambda \left( \delta _{0}\right) =\ln \left[ \delta ^{\left( 3\left(
T+1\right) \right) ^{-1}}\right] \geq \lambda _{1}.  \label{4.46}
\end{equation}%
Then (\ref{4.41})-(\ref{4.46}) imply that%
\begin{equation}
\left. 
\begin{array}{c}
\left\Vert v_{t}\right\Vert _{L_{2}\left( Q_{\varepsilon ,T}\right)
}+\left\Vert \Delta v\right\Vert _{L_{2}\left( Q_{\varepsilon ,T}\right)
}+\left\Vert v\right\Vert _{H^{1,0}\left( Q_{\varepsilon ,T}\right)
}+\left\Vert p\right\Vert _{H^{1,0}\left( Q_{\varepsilon ,T}\right) }\leq \\ 
\leq C_{4}\left( 1+\left\Vert p\right\Vert _{H^{2}\left( Q_{T}\right)
}\right) \delta ^{\rho },\text{ }\forall \delta \in \left( 0,\delta
_{0}\right) ,\text{ }\rho \in \left( 0,1/6\right) .%
\end{array}%
\right.  \label{4.460}
\end{equation}
The rest of the proof of the target H\"{o}lder stability estimates (\ref%
{4.10})-(\ref{4.12}) follows immediately from (\ref{4.7}), (\ref{4.8}), (\ref%
{4.15}), (\ref{4.16}), (\ref{4.460}) and Lemma 3.1.

We now prove uniqueness. Assume that identities (\ref{4.13}) and (\ref{4.130}%
) hold. Then by (\ref{4.9}) and (\ref{4.90}) $\delta =0.$ Hence, (\ref{4.10}%
) and (\ref{4.11}) imply that $u_{1}\left( x,t\right) =u_{2}\left(
x,t\right) $ and $m_{1}\left( x,t\right) =m_{2}\left( x,t\right) $ for $%
\left( x,t\right) \in Q_{\varepsilon ,T}.$ Setting $\varepsilon \rightarrow
0,$ we obtain $u_{1}\left( x,t\right) =u_{2}\left( x,t\right) $ and $%
m_{1}\left( x,t\right) =m_{2}\left( x,t\right) $ for $\left( x,t\right) \in
Q_{T}.$ $\square $

\section{H\"{o}lder Stability and Uniqueness for Problem 2}

\label{sec:5}

Similarly with (\ref{4.0}), for any number $\varepsilon \in \left(
0,T\right) $ define the domain $P_{\varepsilon ,T}$ as:%
\begin{equation}
P_{\varepsilon ,T}=\Omega \times \left( 0,T-\varepsilon \right) \subset
Q_{T}.  \label{5.1}
\end{equation}

\textbf{Theorem 5.1.} \emph{As in Theorem 4.1, let }$%
D_{1},D_{2},D_{3},D_{4}>0$\emph{\ be certain numbers. Let functions }$%
G_{1},G_{2}\in L_{2}\left( Q_{T}\right) $\emph{\ be the right hand sides of
equations (\ref{2.1}), (\ref{2.2})}$.$ \emph{Let functions\ }$F,M,r$ \emph{%
in (\ref{2.1}), (\ref{2.2}) satisfy conditions of Theorem 4.1. Keep
notations (\ref{4.3})-(\ref{4.5}) of Theorem 4.1. Assume that two pairs of
functions }$\left( u_{1},m_{1}\right) $\emph{\ and }$\left(
u_{2},m_{2}\right) $\emph{\ satisfy equations (\ref{2.1})-(\ref{2.3}) with
two pairs of functions }$\left( G_{1,1},G_{2,1}\right) $\emph{\ and }$\left(
G_{1,2},G_{2,2}\right) $\ $\emph{respectively}$\emph{\ and are such that }%
\begin{equation}
\left( u_{1},m_{1}\right) ,\left( u_{2},m_{2}\right) \in K_{3}\left(
D_{3}\right) \times K_{4}\left( D_{4}\right) .  \label{5.2}
\end{equation}%
\emph{Assume that these two pairs of functions }$\left( u_{1},m_{1}\right) $%
\emph{\ and }$\left( u_{2},m_{2}\right) $ \emph{have the following initial
conditions (see (\ref{2.5})):}%
\begin{equation}
u_{1}\left( x,0\right) =u_{0}^{\left( 1\right) }\left( x\right) ,m_{1}\left(
x,0\right) =m_{0}^{\left( 1\right) }\left( x\right) ,\text{ }x\in \Omega .
\label{5.3}
\end{equation}%
\begin{equation}
\text{ }u_{2}\left( x,0\right) =u_{0}^{\left( 2\right) }\left( x\right) ,%
\text{ }m_{2}\left( x,0\right) =m_{0}^{\left( 2\right) }\left( x\right) ,%
\text{ }x\in \Omega .  \label{5.4}
\end{equation}%
\emph{Let the number }$\varepsilon \in \left( 0,T\right) $\emph{\ and let }$%
P_{\varepsilon ,T}$\emph{\ be the domain defined in (\ref{5.1}). Then there
exists a number }$C_{5}=C_{5}\left( \beta ,D,T,\Omega ,\varepsilon \right)
>0 $\emph{\ and a sufficiently small number }$\delta _{0}=$\emph{\ }$\delta
_{0}\left( \beta ,D,T,\Omega ,\varepsilon \right) \in \left( 0,1\right) $ 
\emph{depending only on listed parameters such that if }$\delta \in \left(
0,\delta _{0}\right) $, \emph{\ }%
\begin{equation}
\left\Vert u_{0}^{\left( 1\right) }-u_{0}^{\left( 2\right) }\right\Vert
_{H^{1}\left( \Omega \right) },\left\Vert m_{0}^{\left( 1\right)
}-m_{0}^{\left( 2\right) }\right\Vert _{L_{2}\left( \Omega \right) }\leq
\delta ,  \label{5.5}
\end{equation}%
\emph{and also if inequalities (\ref{4.90}) hold, then there exists a number 
}$\eta =\eta \left( \beta ,D,T,\Omega ,\varepsilon \right) \in \left(
0,1/6\right) $\emph{\ depending only on listed parameters such that the
following H\"{o}lder stability estimate holds:}%
\begin{equation*}
\left. 
\begin{array}{c}
\left\Vert u_{1}-u_{2}\right\Vert _{H^{1,0}\left( P_{\varepsilon ,T}\right)
}+\left\Vert m_{1}-m_{2}\right\Vert _{H^{1,0}\left( P_{\varepsilon
,T}\right) }\leq \\ 
\leq C_{5}\left( 1+\left\Vert u_{1}-u_{2}\right\Vert _{H^{2}\left(
Q_{T}\right) }+\left\Vert m_{1}-m_{2}\right\Vert _{H^{1}\left( Q_{T}\right)
}\right) \delta ^{\eta },\text{ }\forall \delta \in \left( 0,\delta
_{0}\right) .%
\end{array}%
\right.
\end{equation*}

\emph{Furthermore, if in (\ref{5.3}) and (\ref{5.4}) }%
\begin{equation*}
u_{0}^{\left( 1\right) }\left( x\right) \equiv u_{0}^{\left( 2\right)
}\left( x\right) ,\text{ }m_{0}^{\left( 1\right) }\left( x\right) \equiv
m_{0}^{\left( 2\right) }\left( x\right) ,\text{ }x\in \Omega ,
\end{equation*}%
\emph{and if (\ref{4.130}) holds as well, then }$u_{1}\left( x,t\right)
\equiv u_{2}\left( x,t\right) $\emph{\ and }$m_{1}\left( x,t\right) \equiv
m_{2}\left( x,t\right) $\emph{\ in }$Q_{T},$\emph{\ i.e. Problem 2 has at
most one solution }$\left( u,m\right) \in K_{3}\left( D_{3}\right) \times
K_{4}\left( D_{4}\right) $\emph{.}

\textbf{Proof}. Similarly with the proof of Theorem 4.1, in this proof $%
\widetilde{C}_{5}=\widetilde{C}_{5}\left( \beta ,D,T,\Omega \right) >0$
denotes different numbers depending only on $\beta ,D,T,\Omega $ and $%
C_{5}=C_{5}\left( \beta ,D,T,\Omega ,\varepsilon \right) >0$ denotes
different numbers depending not only on parameters $\beta ,D,T,\Omega $ but
on $\varepsilon $ as well. We will choose below the number $c=c\left(
T\right) >2$ in (\ref{3.4}) as:%
\begin{equation}
c=c\left( T\right) =2+\sqrt{T+\frac{1}{4}}.  \label{5.60}
\end{equation}%
We introduce the number $\xi =\xi \left( T\right) ,$ 
\begin{equation}
\xi =\xi \left( T\right) =\frac{T+c}{c^{2}}=\frac{T+2+\sqrt{T+1/4}}{\left( 2+%
\sqrt{T+1/4}\right) ^{2}}\in \left( 0,1\right) .  \label{5.62}
\end{equation}%
The reason of the choice of (\ref{5.60}), (\ref{5.62}) is explained in this
proof below.

Keep notations (\ref{4.15}) and (\ref{4.160}) and replace (\ref{4.16}) with%
\begin{equation}
v\left( x,0\right) =v_{0}\left( x\right) =u_{0}^{\left( 1\right) }\left(
x\right) -u_{0}^{\left( 2\right) }\left( x\right) ,\text{ }x\in \Omega ,
\label{5.7}
\end{equation}%
\begin{equation}
p\left( x,0\right) =p_{0}\left( x\right) =m_{0}^{\left( 1\right) }\left(
x\right) -m_{0}^{\left( 2\right) }\left( x\right) ,\text{ }x\in \Omega .
\label{5.70}
\end{equation}

Similarly with (\ref{4.19})-(\ref{4.22}) we obtain two inequalities:%
\begin{equation}
\left\vert v_{t}+\beta \Delta v\right\vert \left( x,t\right) \leq \widetilde{%
C}_{5}\left( \left\vert \nabla v\right\vert +\dint\limits_{\Omega
}\left\vert p\left( y,t\right) \right\vert dy+\left\vert p\right\vert
+\left\vert \widetilde{G}_{1}\right\vert \right) \left( x,t\right) ,\text{ }%
\left( x,t\right) \in Q_{T},  \label{5.8}
\end{equation}%
\begin{equation}
\left. 
\begin{array}{c}
\left\vert p_{t}-\beta \Delta p+r^{2}\left( x\right) m_{1}\Delta
v\right\vert \left( x,t\right) \leq \\ 
\leq \widetilde{C}_{5}\left( \left\vert \nabla p\right\vert +\left\vert
p\right\vert +\left\vert \nabla v\right\vert +\left\vert \widetilde{G}%
_{2}\right\vert \right) \left( x,t\right) ,\text{ }\left( x,t\right) \in
Q_{T},%
\end{array}%
\right.  \label{5.9}
\end{equation}%
\begin{equation}
\partial _{\nu }v\mid _{S_{T}}=0,\text{ }\partial _{\nu }p\mid _{S_{T}}=0,
\label{5.10}
\end{equation}%
\begin{equation}
v\left( x,0\right) =v_{0}\left( x\right) ,\text{ }p\left( x,0\right)
=p_{0}\left( x\right) .  \label{5.11}
\end{equation}

Note that the difference between (\ref{5.9}) and (\ref{4.20}) is that the
term with $\Delta v$ is in the left hand side of (\ref{5.9}) rather than
being in the right hand side of (\ref{4.20}). This is because the
quasi-Carleman estimate of Theorem 3.4, being applied to the left hand side
of (\ref{5.9}), helps us to handle this.

Below the parameter $c=c\left( T\right) $ in the Carleman Weight Function $%
\varphi _{\lambda }\left( t\right) $ in (\ref{3.4}) is as in (\ref{5.60}).
Applying to (\ref{5.8}) and (\ref{5.9}) the same procedure as the one in the
proof of Theorem 4.1, we obtain the following analogs of (\ref{4.23}) and (%
\ref{4.24}):%
\begin{equation}
\left. 
\begin{array}{c}
\dint\limits_{Q_{T}}\left( v_{t}+\beta \Delta v\right) ^{2}\varphi _{\lambda
}^{2}dxdt\leq \widetilde{C}_{5}\dint\limits_{Q_{T}}\left( \nabla v\right)
^{2}\varphi _{\lambda }^{2}dxdt+ \\ 
+\widetilde{C}_{5}\dint\limits_{Q_{T}}\left( p^{2}+\dint\limits_{\Omega
}p^{2}\left( y,t\right) dy+\widetilde{G}_{1}^{2}\right) \varphi _{\lambda
}^{2}dxdt,%
\end{array}%
\right.  \label{5.12}
\end{equation}%
\begin{equation}
\left. 
\begin{array}{c}
\dint\limits_{Q_{T}}\left( p_{t}-\beta \Delta p+r^{2}\left( x\right)
m_{1}\Delta v\right) ^{2}\varphi _{\lambda }^{2}dxdt\leq \widetilde{C}%
_{5}\dint\limits_{Q_{T}}\left( \left( \nabla p\right) ^{2}+p^{2}\right)
\varphi _{\lambda }^{2}dxdt+ \\ 
+\widetilde{C}_{5}\dint\limits_{Q_{T}}\left( \left( \nabla v\right) ^{2}+%
\widetilde{G}_{2}^{2}\right) \varphi _{\lambda }^{2}dxdt.%
\end{array}%
\right.  \label{5.13}
\end{equation}%
Take $\lambda _{0}=\lambda _{0}\left( T\right) $ as in (\ref{3.7}) and apply
Carleman estimate (\ref{3.8}) to the left hand side of (\ref{5.12}). We
obtain%
\begin{equation}
\left. 
\begin{array}{c}
\widetilde{C}_{5}e^{2c^{\lambda }}\dint\limits_{\Omega }\left( \left( \nabla
v\right) ^{2}+v^{2}\right) \left( x,T\right) +\widetilde{C}_{5}\lambda
\left( T+c\right) ^{\lambda -1}e^{2\left( T+c\right) ^{\lambda
}}\dint\limits_{\Omega }v^{2}\left( x,0\right) dx+ \\ 
+\widetilde{C}_{5}\dint\limits_{Q_{T}}\left( \nabla v\right) ^{2}\varphi
_{\lambda }^{2}dxdt+\widetilde{C}_{5}\dint\limits_{Q_{T}}\left(
p^{2}+\dint\limits_{\Omega }p^{2}\left( y,t\right) dy+\widetilde{G}%
_{1}^{2}\right) \varphi _{\lambda }^{2}dxdt\geq \\ 
\geq \sqrt{\lambda }\dint\limits_{Q_{T}}\left( \nabla v\right) ^{2}\varphi
_{\lambda }^{2}dxdt+\lambda ^{2}c^{\lambda
-2}\dint\limits_{Q_{T}}v^{2}\varphi _{\lambda }^{2}dxdt,\text{ }\forall
\lambda \geq \lambda _{0}.%
\end{array}%
\right.  \label{5.14}
\end{equation}%
Choose $\lambda _{2}=\lambda _{2}\left( \beta ,D,T,\Omega \right) \geq
\lambda _{0}>64$ such that 
\begin{equation}
\sqrt{\lambda _{2}}\geq 2\widetilde{C}_{5}.  \label{5.15}
\end{equation}%
Then the term with $\widetilde{C}_{5}\left( \nabla v\right) ^{2}$ in the
left hand side of (\ref{5.14}) is absorbed by the term with $\sqrt{\lambda }%
\left( \nabla v\right) ^{2}$ in the right hand side of (\ref{5.14}). Hence,
we obtain%
\begin{equation}
\left. 
\begin{array}{c}
\widetilde{C}_{5}e^{2c^{\lambda }}\dint\limits_{\Omega }\left( \left( \nabla
v\right) ^{2}+v^{2}\right) \left( x,T\right) +\widetilde{C}_{5}\lambda
\left( T+c\right) ^{\lambda -1}e^{2\left( T+c\right) ^{\lambda
}}\dint\limits_{\Omega }v^{2}\left( x,0\right) dx+ \\ 
+\widetilde{C}_{5}\dint\limits_{Q_{T}}\left( p^{2}+\dint\limits_{\Omega
}p^{2}\left( y,t\right) dy+\widetilde{G}_{1}^{2}\right) \varphi _{\lambda
}^{2}dxdt\geq \\ 
\geq \sqrt{\lambda }\dint\limits_{Q_{T}}\left( \nabla v\right) ^{2}\varphi
_{\lambda }^{2}dxdt+\lambda ^{2}c^{\lambda
-2}\dint\limits_{Q_{T}}v^{2}\varphi _{\lambda }^{2}dxdt,\text{ }\forall
\lambda \geq \lambda _{2}.%
\end{array}%
\right.  \label{5.16}
\end{equation}%
Since 
\begin{equation*}
\dint\limits_{Q_{T}}\left( \dint\limits_{\Omega }p^{2}\left( y,t\right)
dy\right) \varphi _{\lambda }^{2}dxdt\leq \widetilde{C}_{5}\dint%
\limits_{Q_{T}}p^{2}\varphi _{\lambda }^{2}dxdt,
\end{equation*}%
then (\ref{5.16}) implies%
\begin{equation}
\left. 
\begin{array}{c}
\widetilde{C}_{5}e^{2c^{\lambda }}\dint\limits_{\Omega }\left( \left( \nabla
v\right) ^{2}+v^{2}\right) \left( x,T\right) +\widetilde{C}_{5}\lambda
\left( T+c\right) ^{\lambda -1}e^{2\left( T+c\right) ^{\lambda
}}\dint\limits_{\Omega }v^{2}\left( x,0\right) dx+ \\ 
+\widetilde{C}_{5}\dint\limits_{Q_{T}}\left( p^{2}+\widetilde{G}%
_{1}^{2}\right) \varphi _{\lambda }^{2}dxdt\geq \\ 
\geq \sqrt{\lambda }\dint\limits_{Q_{T}}\left( \nabla v\right) ^{2}\varphi
_{\lambda }^{2}dxdt+\lambda ^{2}c^{\lambda
-2}\dint\limits_{Q_{T}}v^{2}\varphi _{\lambda }^{2}dxdt,\text{ }\forall
\lambda \geq \lambda _{2}.%
\end{array}%
\right.  \label{5.17}
\end{equation}

Now the question is on how to estimate the integral containing $p^{2}$ in
the left hand side of (\ref{5.17}). We will do this via working with (\ref%
{5.13}).

Using Carleman estimate (\ref{3.9}), we estimate the left hand side of (\ref%
{5.13}) as:%
\begin{equation*}
\dint\limits_{Q_{T}}\left( p_{t}-\beta \Delta p+r^{2}\left( x\right)
m_{1}\Delta v\right) ^{2}\varphi _{\lambda }^{2}dxdt\geq
\end{equation*}%
\begin{equation*}
\geq \widetilde{C}_{5}\lambda c^{\lambda -1}\dint\limits_{Q_{T}}\left(
\nabla p\right) ^{2}\varphi _{\lambda }^{2}dxdt+\widetilde{C}_{5}\lambda
^{2}c^{2\lambda -2}\dint\limits_{Q_{T}}p^{2}\varphi _{\lambda }^{2}dxdt-
\end{equation*}%
\begin{equation*}
-\widetilde{C}_{5}\lambda \left( T+c\right) ^{\lambda
-1}\dint\limits_{Q_{T}}\left( \nabla v\right) ^{2}\varphi _{\lambda
}^{2}dxdt-
\end{equation*}%
\begin{equation*}
-\lambda \left( T+c\right) ^{\lambda -1}e^{2\left( T+c\right) ^{\lambda
}}\dint\limits_{\Omega }p^{2}\left( x,0\right) dx,\text{ }\forall \lambda
\geq \lambda _{0}.
\end{equation*}%
Comparing this with (\ref{5.13}), we obtain%
\begin{equation}
\left. 
\begin{array}{c}
\lambda \left( T+c\right) ^{\lambda -1}e^{2\left( T+c\right) ^{\lambda
}}\dint\limits_{\Omega }p^{2}\left( x,0\right) dx+ \\ 
+\widetilde{C}_{5}\dint\limits_{Q_{T}}\left( \left( \nabla p\right)
^{2}+p^{2}\right) \varphi _{\lambda }^{2}dxdt+ \\ 
+\widetilde{C}_{5}\dint\limits_{Q_{T}}\left( \left( \nabla v\right) ^{2}+%
\widetilde{G}_{2}^{2}\right) \varphi _{\lambda }^{2}dxdt+\widetilde{C}%
_{5}\lambda \left( T+c\right) ^{\lambda -1}\dint\limits_{Q_{T}}\left( \nabla
v\right) ^{2}\varphi _{\lambda }^{2}dxdt\geq \\ 
\geq \lambda c^{\lambda -1}\dint\limits_{Q_{T}}\left( \nabla p\right)
^{2}\varphi _{\lambda }^{2}dxdt+\lambda ^{2}c^{2\lambda
-2}\dint\limits_{Q_{T}}p^{2}\varphi _{\lambda }^{2}dxdt,\text{ }\forall
\lambda \geq \lambda _{0}.%
\end{array}%
\right.  \label{5.18}
\end{equation}%
Since $\lambda _{2}\geq \lambda _{0}>64$ and $c>2$ by (\ref{5.60}), then (%
\ref{5.15}) implies%
\begin{equation*}
\lambda \left( T+c\right) ^{\lambda -1}>2\widetilde{C}_{5}\text{ and }%
\lambda c^{\lambda -1}>2\widetilde{C}_{5},\text{ }\forall \lambda \geq
\lambda _{2}.
\end{equation*}%
Hence, (\ref{5.18}) can be rewritten as: 
\begin{equation}
\left. 
\begin{array}{c}
\widetilde{C}_{5}\lambda \left( T+c\right) ^{\lambda -1}e^{2\left(
T+c\right) ^{\lambda }}\dint\limits_{\Omega }p^{2}\left( x,0\right) dx+%
\widetilde{C}_{5}\dint\limits_{Q_{T}}\widetilde{G}_{2}^{2}\varphi _{\lambda
}^{2}dxdt+ \\ 
+\widetilde{C}_{5}\lambda \left( T+c\right) ^{\lambda
-1}\dint\limits_{Q_{T}}\left( \nabla v\right) ^{2}\varphi _{\lambda
}^{2}dxdt\geq \\ 
\geq \lambda c^{\lambda -1}\dint\limits_{Q_{T}}\left( \nabla p\right)
^{2}\varphi _{\lambda }^{2}dxdt+\lambda ^{2}c^{2\lambda
-2}\dint\limits_{Q_{T}}p^{2}\varphi _{\lambda }^{2}dxdt,\text{ }\forall
\lambda \geq \lambda _{2}.%
\end{array}%
\right.  \label{5.19}
\end{equation}%
In particular, (\ref{5.19}) gives us the following desired estimate for the
integral containing $p^{2}:$%
\begin{equation}
\left. 
\begin{array}{c}
\dint\limits_{Q_{T}}p^{2}\varphi _{\lambda }^{2}dxdt\leq \widetilde{C}%
_{5}\xi ^{\lambda -1}\dint\limits_{Q_{T}}\left( \nabla v\right) ^{2}\varphi
_{\lambda }^{2}dxdt+ \\ 
+\widetilde{C}_{5}\left( T+c\right) ^{\lambda -1}e^{2\left( T+c\right)
^{\lambda }}\dint\limits_{\Omega }p^{2}\left( x,0\right) dx+\widetilde{C}%
_{5}\dint\limits_{Q_{T}}\widetilde{G}_{2}^{2}\varphi _{\lambda }^{2}dxdt,%
\text{ }\forall \lambda \geq \lambda _{2},%
\end{array}%
\right.  \label{5.21}
\end{equation}%
where the number $\xi =\xi \left( T\right) \in \left( 0,1\right) $ is
defined in (\ref{5.62}), and this is the reason for our above choice of (\ref%
{5.60}), (\ref{5.62}). Hence,%
\begin{equation}
\left. 
\begin{array}{c}
-\widetilde{C}_{5}\dint\limits_{Q_{T}}p^{2}\varphi _{\lambda }^{2}dxdt\geq
-\xi ^{\lambda -1}\dint\limits_{Q_{T}}\left( \nabla v\right) ^{2}\varphi
_{\lambda }^{2}dxdt- \\ 
-\widetilde{C}_{5}\left( T+c\right) ^{\lambda -1}e^{2\left( T+c\right)
^{\lambda }}\dint\limits_{\Omega }p^{2}\left( x,0\right) dx+\widetilde{C}%
_{5}\dint\limits_{Q_{T}}\widetilde{G}_{2}^{2}\varphi _{\lambda }^{2}dxdt,%
\text{ }\forall \lambda \geq \lambda _{2}.%
\end{array}%
\right.  \label{5.22}
\end{equation}%
Choose $\lambda _{3}=\lambda _{3}\left( \beta ,D,T,\Omega \right) \geq
\lambda _{2}$ such that%
\begin{equation*}
\xi ^{\lambda -1}\leq \frac{\sqrt{\lambda }}{2},\text{ }\forall \lambda \geq
\lambda _{3}.
\end{equation*}%
Then substitute (\ref{5.22}) in (\ref{5.17}). We obtain%
\begin{equation}
\left. 
\begin{array}{c}
\dint\limits_{Q_{T}}\left( \left( \nabla v\right) ^{2}+v^{2}\right) \varphi
_{\lambda }^{2}dxdt\leq \\ 
\leq \widetilde{C}_{5}e^{2c^{\lambda }}\dint\limits_{\Omega }\left( \left(
\nabla v\right) ^{2}+v^{2}\right) \left( x,T\right) +\widetilde{C}%
_{5}e^{3\left( T+c\right) ^{\lambda }}\dint\limits_{\Omega }v^{2}\left(
x,0\right) dx+ \\ 
+\widetilde{C}_{5}e^{3\left( T+c\right) ^{\lambda }}\dint\limits_{\Omega
}p^{2}\left( x,0\right) dx+\widetilde{C}_{5}\dint\limits_{Q_{T}}\left( 
\widetilde{G}_{1}^{2}+\widetilde{G}_{2}^{2}\right) \varphi _{\lambda
}^{2}dxdt,\text{ }\forall \lambda \geq \lambda _{3}.%
\end{array}%
\right.  \label{5.23}
\end{equation}%
Hence, the integral with $\left( \nabla v\right) ^{2}$ in the left hand side
of (\ref{5.19}) can be estimated via the right hand side of (\ref{5.23}).
This means, in turn that integrals in the right hand side of (\ref{5.19})
can also be estimated via the right hand side of (\ref{5.23}). Thus, we
obtain 
\begin{equation}
\left. 
\begin{array}{c}
\dint\limits_{Q_{T}}\left( \left( \nabla p\right) ^{2}+p^{2}\right) \varphi
_{\lambda }^{2}dxdt\leq \\ 
\leq \widetilde{C}_{5}e^{2c^{\lambda }}\dint\limits_{\Omega }\left( \left(
\nabla v\right) ^{2}+v^{2}\right) \left( x,T\right) +\widetilde{C}%
_{5}\lambda \left( T+c\right) ^{\lambda -1}e^{2\left( T+c\right) ^{\lambda
}}\dint\limits_{\Omega }v^{2}\left( x,0\right) dx+ \\ 
+\widetilde{C}_{5}\left( T+c\right) ^{\lambda -1}e^{2\left( T+c\right)
^{\lambda }}\dint\limits_{\Omega }p^{2}\left( x,0\right) dx+\widetilde{C}%
_{5}\dint\limits_{Q_{T}}\left( \widetilde{G}_{1}^{2}+\widetilde{G}%
_{2}^{2}\right) \varphi _{\lambda }^{2}dxdt,\text{ }\forall \lambda \geq
\lambda _{3}.%
\end{array}%
\right.  \label{5.24}
\end{equation}

We recall now the domain $P_{\varepsilon ,T}\subset Q_{T}$ defined in (\ref%
{5.1}). Hence, using (\ref{4.90}), (\ref{5.5}), (\ref{5.7}), (\ref{5.70}), (%
\ref{5.23}), (\ref{5.24}) and the trace theorem, we obtain%
\begin{equation}
\left. 
\begin{array}{c}
\dint\limits_{P_{\varepsilon ,T}}\left( \left( \nabla v\right)
^{2}+v^{2}+\left( \nabla p\right) ^{2}+p^{2}\right) \varphi _{\lambda
}^{2}dxdt\leq \\ 
\leq \widetilde{C}_{5}e^{3\left( T+c\right) ^{\lambda }}\delta ^{2}+%
\widetilde{C}_{5}e^{2c^{\lambda }}\left( \left\Vert v\right\Vert
_{H^{2}\left( Q_{T}\right) }^{2}+\left\Vert p\right\Vert _{H^{1}\left(
Q_{T}\right) }^{2}\right) ,\text{ }\forall \lambda \geq \lambda _{3}.%
\end{array}%
\right.  \label{5.25}
\end{equation}%
By (\ref{3.4}) and (\ref{5.1}) 
\begin{equation*}
\min_{\overline{P}_{\varepsilon ,T}}\varphi _{\lambda }\left( t\right)
=e^{\left( c+\varepsilon \right) ^{\lambda }}.
\end{equation*}%
Hence, using (\ref{5.25}), we obtain%
\begin{equation}
\left. 
\begin{array}{c}
\left\Vert v\right\Vert _{H^{1,0}\left( P_{\varepsilon ,T}\right)
}^{2}+\left\Vert p\right\Vert _{H^{1,0}\left( P_{\varepsilon ,T}\right)
}^{2}\leq C_{5}e^{3\left( T+c\right) ^{\lambda }}\delta ^{2}+ \\ 
+C_{5}\exp \left[ -2\left( c+\varepsilon \right) ^{\lambda }\left( 1-\left(
c/\left( c+\varepsilon \right) \right) ^{\lambda }\right) \right] \left(
\left\Vert v\right\Vert _{H^{2}\left( Q_{T}\right) }^{2}+\left\Vert
p\right\Vert _{H^{1}\left( Q_{T}\right) }^{2}\right) ,\text{ } \\ 
\forall \lambda \geq \lambda _{3}.%
\end{array}%
\right.  \label{5.26}
\end{equation}%
Choose $\lambda _{4}=\lambda _{4}\left( \beta ,D,T,\Omega ,\varepsilon
\right) \geq \lambda _{3}$ such that%
\begin{equation*}
\left( \frac{c}{c+\varepsilon }\right) ^{\lambda }\leq \frac{1}{2},\text{ }%
\forall \lambda \geq \lambda _{4}.
\end{equation*}%
Then (\ref{5.26}) implies%
\begin{equation}
\left. 
\begin{array}{c}
\left\Vert v\right\Vert _{H^{1,0}\left( P_{\varepsilon ,T}\right)
}^{2}+\left\Vert p\right\Vert _{H^{1,0}\left( P_{\varepsilon ,T}\right)
}^{2}\leq \\ 
\leq C_{5}e^{3\left( T+c\right) ^{\lambda }}\delta ^{2}+C_{5}e^{-\left(
c+\varepsilon \right) ^{\lambda }}\left( \left\Vert v\right\Vert
_{H^{2}\left( Q_{T}\right) }^{2}+\left\Vert p\right\Vert _{H^{1}\left(
Q_{T}\right) }^{2}\right) ,\text{ }\forall \lambda \geq \lambda _{4}.%
\end{array}%
\right.  \label{5.27}
\end{equation}

Let $\lambda =\lambda \left( \delta \right) $ be such that 
\begin{equation}
e^{3\left( T+c\right) ^{\lambda }}\delta ^{2}=\delta .  \label{5.28}
\end{equation}%
Hence, 
\begin{equation}
\lambda \left( \delta \right) =\ln \left[ \delta ^{-\left( 3\left(
T+c\right) ^{-1}\right) }\right] .  \label{5.29}
\end{equation}%
Then 
\begin{equation}
e^{-\left( c+\varepsilon \right) ^{\lambda }}=\delta ^{\eta },\text{ }2\eta =%
\frac{c+\varepsilon }{3\left( T+c\right) }<\frac{1}{3}.  \label{5.30}
\end{equation}%
Choose $\delta _{0}=\delta _{0}\left( \beta ,D,T,\Omega ,\varepsilon \right)
\in \left( 0,1\right) $ so small that 
\begin{equation}
\ln \left[ \delta _{0}^{-\left( 3\left( T+c\right) ^{-1}\right) }\right]
\geq \lambda _{4}.  \label{5.31}
\end{equation}%
Hence, using (\ref{5.27})-(\ref{5.31}), we obtain%
\begin{equation*}
\left. 
\begin{array}{c}
\left\Vert v\right\Vert _{H^{1,0}\left( P_{\varepsilon ,T}\right)
}+\left\Vert p\right\Vert _{H^{1,0}\left( P_{\varepsilon ,T}\right) }\leq \\ 
\leq C_{5}\left( 1+\left\Vert v\right\Vert _{H^{2}\left( Q_{T}\right)
}+\left\Vert p\right\Vert _{H^{1}\left( Q_{T}\right) }\right) \delta ^{\eta
},\text{ }\forall \delta \in \left( 0,\delta _{0}\right) .%
\end{array}%
\right.
\end{equation*}%
The rest of the proof is the same as the part of the proof of Theorem 4.1
after (\ref{4.460}). \ $\square $

\end{document}